\documentclass[12pt]{amsart}

\usepackage{amsmath}
\usepackage{amssymb}
\usepackage{latexsym}
\newtheorem{theorem}{Theorem}[section]
\newtheorem{cor}[theorem]{Corollary}
\newtheorem{lemma}[theorem]{Lemma}

\newenvironment{proof*}{\vskip 2mm\noindent {}}{\hfill $\Box$ \vskip 2mm}
\numberwithin{equation}{section}
\newcommand{\C}{{\mathbb{C}}}

\newcommand{\R}{{\mathbb{R}}}

\newcommand{\D}{{\mathbb{D}}}

\newcommand{\eps}{\varepsilon}

\begin{document}

\title{An example of limit of Lempert functions}
                      %\vskip1cm

\author{Pascal J. Thomas}


\begin{abstract} 
The Lempert function for a set of poles in a domain of $\mathbb C^n$ at a point $z$ is obtained by
taking a certain infimum over all analytic disks going through the poles and the point $z$,
and majorizes the corresponding multi-pole pluricomplex Green function. We investigate
the precise behavior of the Lempert function as a set of three poles in $\mathbb C^2$
coalesces to the origin.
\end{abstract} 
\maketitle

\section{Introduction}
Let $\Omega $ be a domain in $\C^n$, and $a_j \in \Omega$, $j=0,...,N$.
Coman's Lempert function is defined by
   \cite{Coman}, \cite{Lempert}:
\begin{multline}
\ell (z):=\ell_{a_0,\dots,a_N} (z):=
\inf \big\{ \sum^N_{j=0} \log|\zeta_j|:
  \varphi(0)=z, 
\\
\varphi(\zeta_j)=a_j, j=0,...,N
\mbox{  for some } \varphi\in \mathcal {O}(\mathbb D,\Omega) \big\} ,
\end{multline}
where $\mathbb D$ is the unit disc in $\mathbb C$.

For most of this paper, we will consider $\Omega = \D^2$, $|z|:= \max(|z_1|,|z_2|)$,
$a_0=(0, 0)$, $a_1=(\eps_1,0)$, and $a_2=(0,\eps_2)$, where $\eps_j\in \D$, $j=1, 2$. 
We will write 
$\ell_\eps(z)$ for the Lempert function with respect to the three poles $a_0, a_1, a_2$ evaluated
at the point $z=(z_1,z_2)\in \D$. It is clear that the Lempert function approaches $-\infty$ near each of its poles $a_j$. 
When $\eps_j \to 0$, all poles concentrate at the origin of $\mathbb C^2$, and the Lempert function may converge to some limit with its singularities concentrated at the origin.

Our goal is to understand in detail the aspect of this singularity. A first remark is that the Lempert function is always related to the corresponding Green function for the same poles,
$$
g(z) := \sup  \left\lbrace u \in PSH(\Omega, \mathbb R_-) :  u(z) \le \log |z-a_j|+C_j, j=0,...,N \right\rbrace ,
$$
where $PSH(\Omega, \mathbb R_-)$ stands for the set of all negative plurisubharmonic functions in $\Omega$.
The inequality $g(z)\le \ell(z)$ always holds, and it is known that it can be strict 
\cite{CarlWieg}, \cite{TraoTh}, \cite{NikoZwo}. If $\ell$ ever turns out to be plurisubharmonic itself, then it must be equal to $g$ \cite{Coman}.  There are some simple cases where limits of sequences of Lempert functions can be
identified as Lempert functions with multiplicities \cite{TraoTh2}.

In the special case that we are studying, the Green function (denoted by $g_\eps$) is not known either, nor do we know whether it always admits a limit as $\eps_j \to 0$, but it can be compared to the Green functions for the three following systems of points:
$$
S_1:= \left\lbrace a_0, a_1\right\rbrace , \quad S_2:= \left\lbrace a_0, a_2\right\rbrace , \quad
S_3:= \left\lbrace a_0, a_1, a_2, (\eps_1, \eps_2) \right\rbrace .
$$
Those are all product sets, so their Green functions are explicitly known 
\cite{Edigarian} as well as their limits when  $\eps_j \to 0$, which are respectively 
\begin{multline*}
g_1(z) := \max (2 \log |z_1|, \log |z_2|), \quad g_2(z) := \max ( \log |z_1|, 2 \log |z_2|), 
\\
g_3(z) := \max ( 2 \log |z_1|, 2 \log |z_2|).
\end{multline*}

Nguyen Van Trao remarked that it follows from the definition of the Green function that 
\begin{equation}
\label{gineq} 
g_3  \le g_\eps \le \min (g_1, g_2)
\end{equation} 
throughout the bidisk, and therefore when $|z_2| \le |z_1|^2$, $g_\eps (z) = g_3(z) = g_1(z) = 2 \log |z_1|$, and 
when $|z_1| \le |z_2|^2$, $g_\eps (z) = g_3(z) = g_2(z) = 2 \log |z_2|$. Also, for any $z$ in the bidisk,
$$
\liminf_{\eps_1, \eps_2 \to 0} \ell_\eps(z) \ge g_\eps(z) \ge g_3(z) = 2  \log |z|.
$$

We first give our result in a special case where the picture is more complete.
\begin{theorem}
\label{simplecase} 
Suppose $\eps_1=\eps_2=\eps$. Then
\begin{enumerate}
\item
If $z_1=0$ or $z_2=0$ or $z_1+z_2=0$, there exists a constant $C>0$ such that
$$
\limsup_{\eps_1, \eps_2 \to 0} \ell_\eps(z) \le 2  \log |z| +C.
$$
\item
For any $c_0>0$, there exists a constant $C=C(c_0)>0$ such that 
for any $z=(z_1,z_2)$ verifying $c_0 \le \left| \frac{z_1}{z_2} \right| \le c_0^{-1}$, and
$c_0 \le \left| 1+  \frac{z_1}{z_2} \right|$, then
\begin{equation*}
\liminf_{\eps \to 0} \ell_\eps(z) \ge \frac{3}{2}  \log |z| - C
\end{equation*} 
and
\begin{equation*}
\liminf_{\eps \to 0} \ell_\eps(z) \le \frac{3}{2}  \log |z| + C.
\end{equation*} 
\end{enumerate}
\end{theorem}

\begin{cor}
If $g$ is some cluster point of the family $\{g_\eps\}$, then for any $z\in \D^2$, 
$$
g(z) \le \frac{3}{2}  \log |z|,
$$
and $g(z_1,-z_1)=2 \log|z_1|$. 
\end{cor} 
This is an improvement over the inequalities \eqref{gineq} when $|z_2|^2 < |z_1| < |z_2|^{1/2}$. Other estimates,
to appear in  an upcoming paper with Jon Magnusson and Ragnar Sigurdsson, show that these new bounds for
the limits of the Green functions
are actually sharp.

\begin{proof}
The function $g$ must be plurisubharmonic, so subharmonic on any complex line going through the origin.
It is also negative everywhere.
The inequality that its restriction on a complex line
satisfies near the origin implies that it is bounded above by the corresponding
one-variable Green function, which gives us the required upper bound. The equality 
in the special case follows from \eqref{gineq}.
\end{proof} 

Theorem \ref{simplecase} will follow from the more detailed result below.

\begin{theorem}
\label{complete} 
\begin{enumerate}
\item
If $z_1=0$ or $z_2=0$ or if $z_1\neq0, z_2\neq 0$ and $\lim_{\eps_1, \eps_2 \to 0} \frac{\eps_1}{\eps_2}=-\frac{z_1}{z_2}$, there exists a constant $C>0$ such that
$$
\limsup_{\eps_1, \eps_2 \to 0} \ell_\eps(z) \le 2  \log |z| +C.
$$
\item
For any $c_0>0$, there exists a constant $C=C(c_0)>0$ such that 
for any $z=(z_1,z_2)$ verifying 
 $c_0 \le \left| \frac{z_1}{z_2} \right| \le c_0^{-1}$, with the $(\eps_1, \eps_2)$ involved in
the limes (inferior and superior) below 
always verifying
\begin{equation}
\label{argcond} 
\left| \frac{z_1}{|z_1|} \frac{|z_2|}{z_2} + \frac{\eps_1}{|\eps_1|} \frac{|\eps_2|}{\eps_2} \right| \ge c_0 ,
\end{equation} 
then 
\begin{equation}
\label{lower} 
\liminf_{\eps_1, \eps_2 \to 0} \ell_\eps(z) \ge \frac{3}{2}\log |z| - C.
\end{equation} 
If, furthermore, 
$$
\lim_{\eps_1, \eps_2 \to 0} \frac{\eps_1^2}{\eps_2} = 
\lim_{\eps_1, \eps_2 \to 0} \frac{\eps_2^2}{\eps_1} = 0,
$$
then
\begin{equation}
\label{upper} 
\limsup_{\eps_1, \eps_2 \to 0} \ell_\eps(z) \le \frac{3}{2}\log |z| + C.
\end{equation} 
\end{enumerate}
\end{theorem}

\section{Upper estimates}
To prove the upper estimates in the above theorem, we shall need to construct appropriate maps
$\varphi$ from the disk to the bidisk. It will be useful to relax a little the condition that $\varphi(\D)\subset \D^2$. 
A more general form of this lemma will appear in  \cite{TraoTh2}.

\begin{lemma}
\label{approx} 
%Suppose $\Omega$ is a bounded hyperconvex domain in $\C^n$, for any $\gamma >0$ let $\Omega^{\gamma}:= %\left\lbrace z \in \C^n : \mbox{dist} (z, \Omega) < \gamma \right\rbrace $.
Suppose that $a_j(\eps) \in \D^2$, $j=0,...,N$ depend on some parameter $\varepsilon \in \C^m$, with $\lim_{\varepsilon \to 0} a_j(\eps) = a_j \in \Omega$, $j=0,...,N$, and that 
 there exists a function $\gamma(\varepsilon)$ such that $\lim_{\varepsilon \to 0}\gamma(\varepsilon)=0$ and $\lambda \in \R$  with the following property: for any $\varepsilon >0$, $z \in \D^2 \setminus \left\lbrace a_1, \dots, a_N \right\rbrace$,  there exists a holomorphic map $\varphi : \D \longrightarrow D(0,1+\gamma(\eps))^2$ and $\zeta_0,\dots,\zeta_N$ such that 
\begin{equation}
\label{map}
\varphi(0)=z, \varphi(\zeta_j)=a_j(\eps), j=0,...,N, \quad \mbox{and } \sum^N_{j=0} \log|\zeta_j| \le \lambda.
\end{equation} 
Then $\limsup_{\varepsilon \to 0} \ell_{a_0(\eps),\dots,a_N(\eps)} (z) \le \lambda$.
\end{lemma}
\begin{proof}
First note that by applying an automorphism $\phi$ of the bidisk exchanging $z$ and $(0,0)$, we have a map $\tilde \varphi:= \phi \circ \varphi $ such that 
$$
\tilde \varphi(0)=(0,0), \varphi(\zeta_j)=\phi(a_j(\eps)), j=0,...,N, \mbox{ and } 
\tilde \varphi : \D \longrightarrow D(0,1+\gamma(\eps))^2,
$$
for another function $\gamma$ with the same property as the original one.

Likewise, to estimate $\ell_{a_1(\eps),\dots,a_N(\eps)} (z)$ it is equivalent to look for maps $\psi$ and points 
$\zeta'_0,\dots,\zeta'_N$ such that $\psi(\D)\subset \D^2$, and 
$$
\psi(0)=(0,0), \varphi(\zeta'_j)=\phi(a_j(\eps)), j=0,...,N, \mbox{and } \limsup_ {\varepsilon \to 0} \sum^N_{j=0} \log|\zeta'_j| \le \lambda
$$
Applying the usual Schwarz lemma to each coordinate of $\tilde \varphi$, we see that $\tilde \varphi(D(0,(1+\gamma(\eps))^{-1}) \subset \D^2$, and therefore 
$$
\psi (\zeta) := \tilde \varphi (\frac{\zeta}{1+\gamma(\eps)}) \quad \mbox{and } 
\zeta'_j :=(1+\gamma(\eps))\zeta_j
$$
will satisfy our requirements.
\end{proof}

Notation: we always write $\phi_a (\zeta) := \frac{a-\zeta}{1-\zeta \bar a}$ for the involutive automorphism of the unit disk exchanging $a$ and $0$. 
Applying an automorphism of $\D$ exchanging $\zeta_0$ and $0$, the original problem of construction of maps as in \eqref{map} is equivalent to finding 
a holomorphic map $\varphi : \D \longrightarrow D(0,1+\gamma(\eps))^2$ and 
new points $\zeta_0,\dots,\zeta_N$ such that 
\begin{equation}
\label{mapzero}
\varphi(0)=(0,0), \varphi(\zeta_j)=a_j(\eps), j=1,...,N, \quad \mbox{and } 
|\zeta_0| + \sum^N_{j=1} \log|\phi_{\zeta_0}(\zeta_j)| \le \lambda.
\end{equation} 

\begin{proof*}{\it Proof of Theorem  \ref{complete}, part (1).}

Assume $z_2=0$. The case $z_1=0$ would be treated in the same way.

By the formulation \eqref{mapzero} of our problem, we need to construct a map $\varphi$, with an appropriate control of the image of $\D$, such that 
\begin{equation}
\label{probz20}
\left\lbrace 
\begin{array}{rcl}
\varphi(0)&=&(0,0)\\
\varphi(\zeta_1)&=&(\eps_1,0)\\
\varphi(\zeta_2)&=&(0,\eps_2)\\
\varphi(\zeta_0)&=&(z_1,0)
\end{array}
\right. 
\end{equation} 
First we choose  $\zeta_1:=\eps_1$, $\zeta_0:=z_1$, and $\zeta_2$ close to $1$, to be specified later. An approximate solution of this interpolation problem, will be given by the following map from $\D$ to $\D^2$:
$$
\varphi^0 (\zeta) := \left(  \zeta \phi_{\zeta_2}(\zeta), \zeta 
\frac{\phi_{\eps_1}(\zeta)}{\phi_{\eps_1}(1)} 
\frac{\phi_{z_1}(\zeta)}{\phi_{z_1}(1)} 
\eps_2 \right) .
$$
The errors with respect to the requirements in \eqref{probz20} are now given by
\begin{equation}
\label{errz20}
\left\lbrace 
\begin{array}{rcll}
\varphi^0(0)-(0,0)&=& (0,0)&\\
\varphi^0(\eps_1)-(\eps_1,0)&= & \left( \eps_1( \phi_{\zeta_2}(\eps_1)-1),0 \right) &=: (E_1,0)\\
\varphi^0(\zeta_2)-(0,\eps_2)&=&  \left( 0,\eps_2 (\zeta_2 
\frac{\phi_{\eps_1}(\zeta_2)}{\phi_{\eps_1}(1)} 
\frac{\phi_{z_1}(\zeta_2)}{\phi_{z_1}(1)} -1) \right) &=:(0,E_2)\\
\varphi^0(z_1)-(z_1,0)&=&  \left( z_1 (\phi_{\zeta_2}(z_1)-1),0 \right)  &=: (E_3,0)
\end{array}
\right. 
\end{equation} 
A computation shows that 
$$
\left| \phi_{\zeta_2}(\zeta)-1 \right| = \left|  \frac{(1-\zeta_2) + \zeta (1-\bar \zeta_2)}{1-\zeta \bar \zeta_2}\right| 
\le \left( \frac{1+|\zeta|}{1-|\zeta|} \right) |1-\zeta_2|,
$$
so that for $|z_1| \le \frac12$, $\eps_1 \le \frac12 |z_1|$, which we may assume, 
\begin{equation}
\label{err13} 
|E_1| \le  \frac53 |1-\zeta_2| |\eps_1|, |E_3| \le 3 |1-\zeta_2| |z_1|.
\end{equation} 
In the same way, one sees that under the above assumption, $|E_2|\preceq |1-\zeta_2| |\eps_2|$.

To get a map satisfying \eqref{probz20}, we subtract from $\varphi^0$ a 
correcting term $ \varphi^1 (\zeta)= (\varphi^1_1 (\zeta), \varphi^1_2 (\zeta))$ obtained by Lagrange interpolation.
More precisely, we choose 
$$
\varphi^1_2 (\zeta)= \frac{\zeta}{\zeta_2}  
\frac{\phi_{\eps_1}(\zeta)}{\phi_{\eps_1}(\zeta_2)} 
\frac{\phi_{z_1}(\zeta)}{\phi_{z_1}(\zeta_2)} E_2,
$$
and $ \varphi^1_1 (\zeta) =  \zeta \phi_{\zeta_2}(\zeta) h(\zeta)$, where $h$ must satisfy
$$
h(\eps_1) = \frac{E_1}{\eps_1 \phi_{\zeta_2}(\eps_1)} =: E'_1, \quad
h(z_1) = \frac{E_3}{z_1 \phi_{\zeta_2}(z_1)} =: E'_3.
$$
We can then set 
$$
h (\zeta) :=  E'_1 \frac{\zeta-z_1}{\eps_1-z_1} 
+ E'_3 \frac{\zeta-\eps_1}{z_1-\eps_1} . 
$$
It is easy to see that $|\varphi^1_2 (\zeta)| \preceq |E_2| \le |\eps| $ for $\zeta_2$ close enough to $1$.
It follows from \eqref{err13} that $|E'_1|, |E'_3| \preceq |1-\zeta_2|$, therefore 
$|\varphi^1_1 (\zeta)| \preceq |1-\zeta_2||z_1|^{-1} \le |\eps| $ for $\zeta_2$ close enough to $1$. So the map $\varphi^0-\varphi^1$ satisfies the hypotheses of Lemma \ref{approx} and since we have 
$$
|\zeta_0| + \sum^2_{j=1} \log|\phi_{\zeta_0}(\zeta_j)| \le 
\log |z_1| + \log |\phi_{z_1}(\eps_1)| \le \log  |z_1| + \log (|z_1| +|\eps_1|),
$$
we may take $\lambda = 2 \log  |z_1| +\eta$, for any $\eta>0$, which implies the conclusion of Theorem  \ref{complete}, part (1), with $C=0$.
\end{proof*}

\begin{proof*}{\it Proof of Theorem  \ref{complete}, part (2),  \eqref{upper}.}

We will need a few notations. We set $\mu:=z_2/z_1$, 
$\sigma:=\eps_2/\eps_1$. Exchanging coordinates if needed, we may assume
$|\mu|\le 1$, and therefore $|z|=|z_1|$. The hypothesis in the theorem is that 
$$
c_0 \le \left|\frac{|\mu|}{\mu}+ \frac{|\sigma|}{\sigma}\right| =
  \left| 1 +\left( \frac{|\mu|}{\mu}\right)^{-1} \frac{|\sigma|}{\sigma}\right|,
$$
so the complex number $\sigma/\mu$ lies outside of certain plane sector containing $-1$, and in particular
$|1+(\sigma/\mu)|\ge c_0$. We choose a complex number $\nu$ such that
$\nu^{2}:= \left( 1+(\sigma/\mu)\right)^{-1}$; this remains bounded. 

We choose a complex number $\zeta_0$ such that $\zeta_0^2=z_1$.  This means that 
$z=(\zeta_0^2, \mu \zeta_0^2)$. 
We also choose, for each $\eps_1$, an $\eps'_1$ such that ${\eps'_1}^2=\eps_1$. 
Now we set
\begin{equation}
\label{defzeta} 
\zeta_1 := \nu \eps'_1 = \left( 1+\frac{\eps_2}{\eps_1 \mu}\right)^{-1/2} \eps_1^{1/2}, \quad
\zeta_2 := -\frac{\sigma}{\mu}\zeta_1 = - \frac{\eps_2}{\eps_1^{1/2}} \left( 1+\frac{\eps_2}{\eps_1 \mu}\right)^{-1/2}.
\end{equation}
We will follow the pattern of the previous proof. In order to apply Lemma \ref{approx}, we need to produce a map 
satisfying
\begin{equation}
\label{prob}
\left\lbrace 
\begin{array}{rcl}
\varphi(0)&=&(0,0)\\
\varphi(\zeta_1)&=&(\eps_1,0)\\
\varphi(\zeta_2)&=&(0,\eps_2)\\
\varphi(\zeta_0)&=&(z_1,z_2)
\end{array}
\right. 
\end{equation} 
We choose $\zeta_0, \zeta_1, \zeta_2$ as above and set 
$$
\varphi^0 (\zeta)= \left( \zeta \frac{\zeta-\zeta_2}{1-\zeta\bar \zeta_2}, \mu \zeta \frac{\zeta-\zeta_1}{1-\zeta\bar \zeta_1}
\right) .
$$
We remark that our $\zeta_j$ have been chosen so that
 \begin{equation*}
\zeta_1(\zeta_1-\zeta_2) = \eps_1, \mu \zeta_2(\zeta_2-\zeta_1)=\eps_2, \mbox{ and }
|\zeta_1 \zeta_2|= |\nu|^2 |\mu|^{-1}|\eps_2| \le c_0^{-2} |\eps_2| .
\end{equation*} 
Elementary computations yield
\begin{equation}
\label{errz}
\left\lbrace 
\begin{array}{rcll}
\varphi^0(0)-(0,0)&=& (0,0)&\\
\varphi^0(\zeta_1)-(\eps_1,0)&= & 
\left(\eps_1 \dfrac{\zeta_1 \bar \zeta_2}{1-\zeta_1 \bar \zeta_2},0 \right) &=: (E_1,0)\\
\varphi^0(\zeta_2)-(0,\eps_2)&=&  \left( 0, 
\eps_2 \dfrac{ \bar \zeta_1 \zeta_2}{1-\bar \zeta_1 \zeta_2}
\right) &=:(0,E_2)\\
\varphi^0(\zeta_0)-z &=&  \left( 
\zeta_0 \dfrac{\zeta_0^2 \bar \zeta_2 -\zeta_2  }{1-\zeta_0 \bar \zeta_2} ,
\mu \zeta_0 \dfrac{\zeta_0^2 \bar \zeta_1 -\zeta_1  }{1-\zeta_0 \bar \zeta_1}
\right)  &=: (E_3, E_4)
\end{array}
\right.  .
\end{equation} 
We construct a correcting term $\varphi^1 $ by Lagrange interpolation:
$$
\begin{array}{rcl}
\varphi^1_1 (\zeta) &=& E_1 \frac{\zeta(\zeta-\zeta_2)(\zeta-\zeta_0)}{\zeta_1(\zeta_1-\zeta_2)(\zeta_1-\zeta_0)}
+ E_3  \frac{\zeta(\zeta-\zeta_1)(\zeta-\zeta_2)}{\zeta_0(\zeta_0-\zeta_1)(\zeta_0-\zeta_2)},\\
\varphi^1_2 (\zeta) &=& E_2 \frac{\zeta(\zeta-\zeta_1)(\zeta-\zeta_0)}{\zeta_2(\zeta_2-\zeta_1)(\zeta_2-\zeta_0)}
+ E_4  \frac{\zeta(\zeta-\zeta_1)(\zeta-\zeta_2)}{\zeta_0(\zeta_0-\zeta_1)(\zeta_0-\zeta_2)}.
\end{array}
$$
The map $\varphi^0-\varphi^1$ will assume the correct values, now we need to see that it sends $\D$ to a neighborhood of the bidisk by estimating the size of $\varphi^1$. 

First note that for $|\eps|$ small enough, 
$|E_3| \le 3 |z_1|^{1/2}|\zeta_2| \le 3 c_0^{-3/2}|z_1|^{1/2}|\eps_1|^{-1/2} |\eps_2|$,
and $|E_4| \le 3 |z_1|^{1/2}|\zeta_1| \le 3 c_0^{-1/2}|z_1|^{1/2} |\eps_1|^{1/2}$,  therefore, using the last hypothesis of the theorem, for $|\eps|$ small enough (depending on $|z_1|$), the second terms in $\varphi^1_1 (\zeta)$ and
$\varphi^1_2 (\zeta)$ can be made arbitrarily small.

On the other hand, for $|\eps|$ small enough, 
\begin{multline*}
\left| \frac{E_1}{\zeta_1(\zeta_1-\zeta_2)}\right| = \left|  \dfrac{\zeta_1 \bar \zeta_2}{1-\zeta_1 \bar \zeta_2}\right|
\le 2 c_0^{-2} |\eps_2|, \\
\left| \frac{E_2}{\zeta_2(\zeta_2-\zeta_1)}\right| = \left|  \dfrac{ \mu \bar \zeta_1\zeta_2}{1- \bar \zeta_1 \zeta_2}\right|
\le 2 c_0^{-3} |\eps_2|.
\end{multline*}
So the first terms in $\varphi^1_1 (\zeta)$ and
$\varphi^1_2 (\zeta)$ are bounded above by $2 c_0^{-3} |\eps_2||z_1|^{-1/2}$, which once again can be made arbitrarily small
for $|\eps|$ small enough. 

Finally, the relevant sum
$$
\log |\zeta_0| + \log |\phi_{\zeta_0}(\zeta_1)| + \log |\phi_{\zeta_0}(\zeta_2)| 
\le \log |\zeta_0| + \log( |\zeta_0| + |\zeta_1|) +  \log( |\zeta_0| + |\zeta_2|) ,
$$
so that $\lambda = 3 \log |\zeta_0| +\eta = \frac{3}{2} \log|z_1| +\eta \le \frac{3}{2} \log|z| +\eta $
may be used to apply Lemma \ref{approx}, for any $\eta>0$.
\end{proof*}

\section{Lower Estimates}

\begin{proof*}{\it Proof of Theorem  \ref{complete}, part (2),  \eqref{lower}.}

We will assume that the conclusion fails, i.e. that for any $C_H>0$ there exist arbitrarily small values of $\eps = (\eps_1,\eps_2)$ such that 
\begin{equation}
\label{contrhyp}
 \ell_\eps(z) \le \frac{3}{2}\log |z| - C_H,
\end{equation} 
which means (we change the value of the constant slightly while keeping the same notation) that there exists a 
holomorphic map $\varphi$ from $\D$ to $\D^2$ and points $\zeta_j \in \D$ satisfying the conditions in \eqref{prob} with 
\begin{equation}
\label{hyplemp} 
\log |\zeta_0| + \log |\phi_{\zeta_0}(\zeta_1)| + \log |\phi_{\zeta_0}(\zeta_2)|  \le  \frac{3}{2}\log |z| - C_H.
\end{equation} 
The interpolation conditions in \eqref{prob} are equivalent to the existence of two holomorphic functions $h_1$, $h_2$ from
$\D$ to itself such that 
\begin{equation*}
\varphi (\zeta) = \left( \zeta \phi_{\zeta_2} (\zeta) h_1(\zeta), \zeta \phi_{\zeta_1} (\zeta) h_2 (\zeta) \right) ,
\end{equation*} 
such that furthermore
\begin{eqnarray}
\label{h11} 
h_1(\zeta_1)&=& \frac{\eps_1}{\zeta_1 \phi_{\zeta_2} (\zeta_1)}=:w_1, \\
\label{h12} 
h_1(\zeta_0)&=& \frac{z_1}{\zeta_0 \phi_{\zeta_2} (\zeta_0)}=:w_2, \\
\label{h21} 
h_2(\zeta_2)&=& \frac{\eps_2}{\zeta_2 \phi_{\zeta_1} (\zeta_2)}=:w_4, \\
\label{h22} 
h_2(\zeta_0)&=& \frac{z_2}{\zeta_0\phi_{\zeta_1} (\zeta_0)}=:w_3.
\end{eqnarray} 

For convenience, we will use the invariant (pseudohyperbolic) distance between points of the unit disk given by 
$$
d_G (a,b) := \left| \phi_a (b) \right| =  \left| \phi_b (a) \right| .
$$
By the invariant Schwarz Lemma, the existence of a holomorphic  function $h_1$ mapping $\D$ to itself and satisfying 
\eqref{h11} and \eqref{h12} is equivalent to 
\begin{equation}
\label{ineqex1} 
\left| w_1 \right| < 1, 
\left| w_2 \right| < 1, 
\mbox{ and }
d_G \left( w_1, w_2 \right) < d_G \left( \zeta_1, \zeta_0 \right) = \left| \phi_{\zeta_1} (\zeta_0) \right|.
\end{equation} 
In the same way, the existence of $h_2$ is equivalent to 
\begin{equation}
\label{ineqex2} 
\left| w_3 \right| < 1, 
\left| w_4 \right| < 1, 
\mbox{ and }
d_G \left( w_3, w_4 \right) < d_G \left( \zeta_2, \zeta_0 \right) = \left| \phi_{\zeta_2} (\zeta_0) \right|.
\end{equation} 

The proof will proceed as follows: one main elementary tool is the fact that if the invariant distance 
of two points in the unit disk is small with respect to 
one of their moduli, then the 
difference of their arguments (or rather, the 
distance between their projections to the unit circle) must be small (Lemma \ref{argdg}).

Then, using \eqref{contrhyp} and the fact that $|w_2|$ and $|w_3|$ are in the unit disk, we will show that 
$|\phi_{\zeta_1} (\zeta_0)|$ and $|\phi_{\zeta_2} (\zeta_0)|$ must both be relatively small. This will imply
that $|w_2|$ and $|w_3|$ 	are also relatively big, and because of \eqref{ineqex1}, \eqref{ineqex2}, so will 
be $|w_1|$ and $|w_4|$, and because $\eps$ is small, this will force $|\phi_{\zeta_2} (\zeta_1)|$ to be
small too, and therefore $|\arg(\zeta_1/\zeta_2)|$. This will allow us to show that $\arg(w_1/w_4)$ is close to 
$\arg(-\eps_1/\eps_2)$, because $\phi_{\zeta_2} (\zeta_1)$ is almost opposite to $\phi_{\zeta_1} (\zeta_2)$. 
On the other hand, use of the triangle inequality on the unit circle will show that 
$\arg(w_1/w_4)=\arg[(w_1/w_2)(w_2/w_3)(w_3/w_4)]$ is close to $\arg(z_1/z_2)$. Hypothesis \eqref{argcond} 
will lead to a contradiction when $C_H$ is big enough.

\begin{lemma}
\label{argdg} 
Suppose that $a, b \in \D$, and that $d_G(a,b)=\delta < |a|$, $\delta\le \frac{1}{2}$. Then 
$$
\left| \arg \left( \frac{a}{b} \right) \right| \le 
\arcsin \left( \frac{\delta (1-|a|^2)}{|a|(1-\delta^2)}\right)
\le \frac{3\delta}{|a|} .
$$
\end{lemma} 
Note that in this result, and the computations that follow, we implicitly only consider arguments in the range $[-\frac{\pi}{2}, \frac{\pi}{2}]$.

\begin{proof}
The set of all points $b \in \D$ such that $d_G(a,b)\le \delta$ is a disk of center 
$\gamma:=\frac{a(1-\delta^2)}{1-|a|^2 \delta^2}$ 
and of radius $\rho:=\frac{\delta (1-|a|^2)}{1-|a|^2 \delta^2}$ \cite[Chapter I, Section 1, p. 3]{Garnet}. 
Elementary geometry shows that the
absolute value of the sine of the angle between $a$ and $b$ is bounded by $\rho/|\gamma|$.
\end{proof}

Our starting remark is that \eqref{hyplemp} can be rewritten as 
\begin{multline}
-\log |w_2|-\log |w_3|=
 \log \left| \frac{\zeta_0 \phi_{\zeta_1}(\zeta_0)}{z_2} \right| +
 \log \left| \frac{\zeta_0 \phi_{\zeta_0}(\zeta_2)}{z_1} \right|
\\
\le \log |\zeta_0| + \frac{3}{2} \log |z| - \log|z_1|  - \log|z_2| - C_H
\\
\le  \log |\zeta_0| - \frac{1}{2} \log |z| - \log c_0 -C_H .
\end{multline} 
Since both terms on the left hand side of the first inequality sign
are positive by \eqref{ineqex1}, \eqref{ineqex2}, they are each bounded by the right hand side, and as a consequence, writing $C'_H:=\exp(C_H)$, 
\begin{eqnarray}
\label{inphi10} 
\left|  \phi_{\zeta_1}(\zeta_0)\right| &< & \frac{|z|^{3/2}}{|z_1| C'_H} \le \frac{|z|^{1/2}}{c_0 C'_H},\\
\label{inphi20} 
\left|  \phi_{\zeta_2}(\zeta_0)\right| &<&  \frac{|z|^{3/2}}{|z_2| C'_H} \le \frac{|z|^{1/2}}{c_0 C'_H}.
\end{eqnarray} 
Furthermore, since the right hand side must also be positive, we have
\begin{equation}
\label{inzeta0} 
|\zeta_0| > c_0 C'_H |z|^{1/2} \ge (c_0 C'_H )^2 \max(\left|  \phi_{\zeta_1}(\zeta_0)\right|, \left|  \phi_{\zeta_2}(\zeta_0)\right|).
\end{equation} 
\vskip.5cm
{\it Lower bounds for $|w_2|$ and $|w_3|$, and $|w_1|$ and $|w_4|$.}

By \eqref{inphi20},
\begin{equation}
\label{lbw2} 
|w_2| =\frac{|z_1|}{|\zeta_0 ||\phi_{\zeta_2} (\zeta_0)|} > \frac{c_0 C'_H |z|^{1/2}}{|\zeta_0 |},
\end{equation} 
and since $d_G(w_2,w_1)<| \phi_{\zeta_1} (\zeta_0)|$, classical facts about the invariant distance 
\cite[Lemma 1.4, p. 4]{Garnet} imply
\begin{equation}
\label{lbw1} 
|w_1| \ge \frac{|w_2|-| \phi_{\zeta_1} (\zeta_0)|}{1-|w_2|| \phi_{\zeta_1} (\zeta_0)|}
\ge
\frac{c_0 C'_H |z|^{1/2}}{|\zeta_0 |}-\frac{|z|^{1/2}}{c_0 C'_H}
\ge
 \frac{c_0 C'_H |z|^{1/2}}{2 |\zeta_0 |}
\end{equation} 
when $C_H$ is big enough (depending on $c_0$). 
The same computations go through for $w_3$ and $w_4$ respectively.

\vskip.5cm
{\it Upper bound for $d_G(\zeta_1, \zeta_2)= |\phi_{\zeta_2} (\zeta_1)|$.}

By definition of $w_1$, we have $\phi_{\zeta_2} (\zeta_1)=\eps_1/(\zeta_1 w_1)$. We estimate the modulus of
$|\zeta_1|$ from below by noting that $\zeta_1$ is close to $\zeta_0$, by \eqref{inphi10} and using \eqref{inzeta0}:
\begin{equation}
\label{lbzeta1} 
|\zeta_1| \ge \frac{|\zeta_0|-| \phi_{\zeta_1} (\zeta_0)|}{1-|\zeta_0|| \phi_{\zeta_1} (\zeta_0)|}
\ge  \frac{1}{2} |\zeta_0| ,
\end{equation} 
 when $C_H$ is big enough (depending on $c_0$). So by \eqref{lbw1}
\begin{equation}
\label{inphi12} 
|\phi_{\zeta_2} (\zeta_1)| \le 4\frac{ |\eps_1|}{c_0 C'_H |z|^{1/2}}.
\end{equation} 
Since $|\eps|$ will tend to $0$ as $z$ remains fixed, this estimate is much stronger than the previous \eqref{inphi10}
and \eqref{inphi20}.
From this, \eqref{lbzeta1} and Lemma \ref{argdg}, we deduce 
\begin{equation}
\label{argzeta12} 
\left|  \arg \left( \frac{\zeta_1}{\zeta_2}\right) \right| \le 24 \frac{|\eps_1|}{(c_0 C'_H)^2 |z|}.
\end{equation} 

\vskip.5cm
{\it Argument estimates for $w_1/w_4$.}

First compute
$$
\frac{w_1}{w_4} = \frac{\eps_1}{\eps_2}\frac{\zeta_2}{\zeta_1} \frac{\phi_{\zeta_1} (\zeta_2)}{\phi_{\zeta_2} (\zeta_1)}
=  - \frac{\eps_1}{\eps_2} \frac{\zeta_2}{\zeta_1} \frac{1-\bar \zeta_2 \zeta_1}{1-\bar \zeta_1 \zeta_2}.
$$
We need to bound $\arg (1-\bar \zeta_2 \zeta_1)$. We use the fact that $\zeta_1$ and  $\zeta_2$ are close to 
each other :
$$
1-\bar \zeta_2 \zeta_1 = (1-|\zeta_1|^2 ) \left( 1
-\frac{(\bar \zeta_2 -\bar \zeta_1) }{1-|\zeta_1|^2}  \zeta_1
\right) ,
$$
and we know that $\frac{a-b}{1-|a|^2} = O(d_G(a,b))$, so that \eqref{inphi12}
 implies that 
\begin{equation*}
\left| \arg (1-\bar \zeta_2 \zeta_1)\right| \le c_1 \frac{ |\eps_1|}{c_0 C'_H |z|^{1/2}},
\end{equation*} 
where $c_1$ is some absolute constant. To avoid problems of definition of arguments, use the notation $w^*:=w/|w|$ 
for any nonzero complex number. We have proved that
\begin{equation}
\label{argest14} 
\left|  \left( \frac{w_1}{w_4}\right)^*-  \left( -\frac{\eps_1}{\eps_2}\right)^* \right| \le c_1 \frac{ |\eps_1|}{c_0 C'_H |z|^{1/2}},
\end{equation} 
where $c_1$ is some absolute constant (not the same as above).

\vskip.5cm
{\it Argument estimates for $w_2/w_3$.}

First note that 
$$
\frac{w_2}{w_3} = \frac{z_1}{z_2} \frac{\phi_{\zeta_1} (\zeta_0)}{\phi_{\zeta_2} (\zeta_0)}
=  \frac{z_1}{z_2} \frac{\phi_{\zeta_1} (\zeta_0)}{\phi_{\zeta_0} (\zeta_1)} 
 \frac{\phi_{\zeta_0} (\zeta_1)}{\phi_{\zeta_0} (\zeta_2)}
 \frac{\phi_{\zeta_0} (\zeta_2)}{\phi_{\zeta_2} (\zeta_0)}.
$$
Since $\phi_{\zeta_0} $ is an automorphism of the unit disk, it preserves the invariant distances, so
$d_G(\phi_{\zeta_0} (\zeta_1),\phi_{\zeta_0} (\zeta_2))= d_G (\zeta_1,\zeta_2)$. To apply Lemma \ref{argdg}, we now need a \emph{lower} bound for $|\phi_{\zeta_0} (\zeta_1)|$, say. But the fact that $|w_3|<1$ (see \eqref{ineqex2}) already 
implies that $|\phi_{\zeta_0} (\zeta_1)| > |z_2| >c_0^{-1}|z|$. Finally, using the Lemma with \eqref{inphi12},
\begin{equation}
\label{argestphi12}
\left| \arg  \left(\frac{\phi_{\zeta_0} (\zeta_1)}{\phi_{\zeta_0} (\zeta_2)} \right) \right|  
%\le \arcsin 
%\left(
%\frac{ 4\frac{ |\eps_1|}{c_0 C'_H |z|^{1/2}} \left( 1- c_0^{-2}|z|^2 \right) }{c_0^{-1}|z|\left( 1-\left(  4\frac{ |\eps_1|}{c_0 C'_H |z|^{1/2}}\right)^2 \right) }
% \right) 
\le
12 \frac{ |\eps_1|}{C'_H |z|^{3/2}},
\end{equation} 
provided that $|\eps|$ is small enough. 

Now computations as in the previous paragraph yield
\begin{equation*}
\frac{\phi_{\zeta_1} (\zeta_0)}{\phi_{\zeta_0} (\zeta_1)} 
 \frac{\phi_{\zeta_0} (\zeta_2)}{\phi_{\zeta_2} (\zeta_0)}
=
\frac{1-\bar \zeta_0 \zeta_1}{1-\bar \zeta_1 \zeta_0} 
\frac{1-\bar \zeta_2 \zeta_0}{1-\bar \zeta_0 \zeta_2}.
\end{equation*} 
Since 
$$
(1-\bar \zeta_0 \zeta_1)(1-\bar \zeta_2 \zeta_0) = (1-\bar \zeta_0 \zeta_1)
(1-(\bar \zeta_1 +\bar \zeta_2-\bar \zeta_1) \zeta_0) 
=
|1-\bar \zeta_0 \zeta_1|^2 \left( 1+ \frac{(\bar \zeta_2-\bar \zeta_1) \zeta_0}{1-\bar \zeta_1 \zeta_0}\right) ,
$$
we can estimate $\arg \left( \frac{1-\bar \zeta_0 \zeta_1}{1-\bar \zeta_1 \zeta_0} 
\frac{1-\bar \zeta_2 \zeta_0}{1-\bar \zeta_0 \zeta_2} \right) $ by a constant multiple of $d_G(\zeta_1,\zeta_2)$, 
and as in the previous paragraph,  
\begin{equation*}
\left|  \arg \left( \frac{\phi_{\zeta_1} (\zeta_0)}{\phi_{\zeta_0} (\zeta_1)} 
 \frac{\phi_{\zeta_0} (\zeta_2)}{\phi_{\zeta_2} (\zeta_0)}
\right) \right| \le c_1 \frac{ |\eps_1|}{c_0 C'_H |z|^{1/2}},
\end{equation*} 
changing the value of $c_1$ as needed.

Using this and \eqref{argestphi12}, the final result of this paragraph is then that 
\begin{equation}
\label{argestwz} 
\left|  \left( \frac{w_2}{w_3}\right)^*-  \left( \frac{z_1}{z_2}\right)^* \right| \le 
c_1 \frac{ |\eps_1|}{c_0 C'_H |z|^{1/2}} + 12 \frac{ |\eps_1|}{C'_H |z|^{3/2}} 
\le 
c_2  \frac{ |\eps_1|}{c_0 C'_H |z|^{3/2}}.
\end{equation} 

\vskip.5cm
{\it Triangle inequality, and contradiction.}

The conditions about invariant distances in \eqref{ineqex1} and \eqref{ineqex2} express the fact that $w_1$ is 
close to $w_2$, and that $w_4$ is close to $w_3$. In fact, applying Lemma \ref{argdg} and the estimates \eqref{lbw2} and 
\eqref{inphi10}, we obtain 
\begin{equation*}
\left| \arg  \left( \frac{w_1}{w_2}\right) \right|  \le \arcsin \left(   \frac{\frac{|z|^{1/2}}{c_0 C'_H} \left( 1- \frac{(c_0 C'_H)^2 |z|}{|\zeta_0 |^2}\right) }{\frac{c_0 C'_H |z|^{1/2}}{|\zeta_0 |} \left( 1- \frac{|z|}{(c_0 C'_H)^2} \right) } \right) 
\le  \frac{3}{(c_0 C'_H)^2},
\end{equation*} 
provided that $C'_H$ is large enough (depending on $c_0$).
As before, the same reasoning applies to $w_4$ and $w_3$, and we finally have 
\begin{equation*}
\left|  \left( \frac{w_1}{w_2}\right)^*- 1  \right| \le  \frac{3}{(c_0 C'_H)^2}, \quad
\left|  \left( \frac{w_3}{w_4}\right)^*- 1  \right| \le  \frac{3}{(c_0 C'_H)^2}.
\end{equation*} 

Now we combine this with \eqref{argestwz} (and the fact that rotations are isometries) to obtain
\begin{multline*}
\left|  \left( \frac{w_1}{w_4}\right)^*-  \left( \frac{z_1}{z_2}\right)^* \right| 
=
\left|  \left( \frac{w_1}{w_2}\right)^*  \left( \frac{w_3}{w_4}\right)^*  \left( \frac{w_2}{w_3}\right)^*
-  \left( \frac{z_1}{z_2}\right)^* \right|
\\
\le 
\left|  \left( \frac{w_1}{w_2}\right)^*- 1  \right| + \left|  \left( \frac{w_3}{w_4}\right)^*- 1  \right| 
+ \left|  \left( \frac{w_2}{w_3}\right)^*-  \left( \frac{z_1}{z_2}\right)^* \right| 
\\
\le
  \frac{6}{(c_0 C'_H)^2} + c_2  \frac{ |\eps_1|}{c_0 C'_H |z|^{3/2}}. 
\end{multline*} 
Combining this with \eqref{argest14}, we see that, for $|\eps|$ small enough (depending on $|z|$),
\begin{equation*}
\left|  \left( \frac{z_1}{z_2}\right)^*-  \left( -\frac{\eps_1}{\eps_2}\right)^* \right| \le 
  \frac{6}{(c_0 C'_H)^2} + c_1 \frac{ |\eps_1|}{c_0 C'_H |z|^{1/2}}  + c_2  \frac{ |\eps_1|}{c_0 C'_H |z|^{3/2}}. 
\end{equation*} 
However, our hypothesis \eqref{argcond} precisely says that the left hand side of this must be greater than $c_0$.
Taking $C'_H$ large enough (depending on $c_0$ only), we may assume $\frac{6}{(c_0 C'_H)^2}  < \frac{c_0}{2}$,
so that by taking $|\eps|$ small enough (depending on $|z|$), we obtain a contradiction, q.e.d.

\end{proof*}

\section{The remaining case}

\begin{proof*}{\it Proof of Theorem  \ref{complete}, part (1), when 
$\lim_{\eps_1, \eps_2 \to 0} \frac{\eps_1}{\eps_2}=-\frac{z_1}{z_2}$.}

We shall reuse the notation $\mu=z_1/z_2$, and reduce ourselves to the case $|\mu|\le1$. The hypothesis implies that 
$\varepsilon_2= -(\mu-\gamma)\varepsilon_1$, where $\gamma=\gamma(\varepsilon)$ tends to $0$ as $\varepsilon$ tends to $0$.

This time, instead of constructing an explicit map satisfying an approximate version of our interpolation problem, we will provide values of $\zeta_0$, $\zeta_1$ and $\zeta_2$ such that the conditions \eqref{ineqex1} and \eqref{ineqex2} are satisfied. From now on we take $\zeta_0=1/2$, $\zeta_1= \zeta_0 + \xi$, $\zeta_2= \zeta_1 + \xi'$, where 
\begin{equation*}
\xi= C_1 z_1, \quad
\xi'= \dfrac{\varepsilon_1}{z_1} \dfrac{\zeta_0}{\zeta_1} \dfrac{1-| \zeta_1 | ^2}{1-\zeta_0 \bar \zeta_1} \xi,
\end{equation*} 
and $C_1=40$. 

{\it Standing assumptions.}

Throughout, we will assume that $|z_1|$ is small enough so that $|\xi|\le 1/4$, so that $1/4 \le \mbox{Re} \zeta_1 \le |\zeta_1| \le 3/4$.  This implies that $1/2 \le  |1-\zeta_0 \bar \zeta_1| \le 1$, and therefore
\begin{equation}
\label{xiprimexi} 
|\xi'|\le 4 |\varepsilon_1| \frac{|\xi|}{|z_1|} \le \frac{1}{2}|\xi|,
\end{equation} 
for $\varepsilon$ small enough (depending on $z_1$). Also note that $|\xi'|\le  |\varepsilon_1|/10$. 
The estimate  \eqref{xiprimexi} implies that $|\xi + \xi'| \ge \frac{1}{2}|\xi|$, and also that 
$$
\frac{1}{2}\le |1-\zeta_0 \bar \zeta_2| \le 1, \quad
\frac{1}{4}|1-\zeta_2 \bar \zeta_1| \le 1.
$$

Now we compute the invariant distances between the $\zeta_j$'s. We have
\begin{multline}
\label{autos} 
\phi_{\zeta_1}(\zeta_0) = \dfrac{\xi}{1-\zeta_0 \bar \zeta_1}, \quad
\phi_{\zeta_2}(\zeta_0) = \dfrac{\xi+\xi'}{1-\zeta_0 \bar \zeta_2}, \\
\phi_{\zeta_2}(\zeta_1) = \dfrac{\xi'}{1-\zeta_1 \bar \zeta_2}, \quad
\phi_{\zeta_1}(\zeta_2) = \dfrac{-\xi'}{1-\zeta_2 \bar \zeta_1}.
\end{multline} 
This implies in particular that 
\begin{equation}
\label{autoest} 
\left| \phi_{\zeta_1}(\zeta_0) \right| \ge |\xi|= C_1 |z_1|, \left| \phi_{\zeta_2}(\zeta_0) \right| \ge |\xi'|\ge \frac{1}{2}|\xi|
=\frac{C_1}{2}  |z_1|.
\end{equation} 

Then 
\begin{multline*}
\log |\zeta_0| + \log |\phi_{\zeta_0}(\zeta_1)| + \log |\phi_{\zeta_0}(\zeta_2)| \\
= -\log  2 + \log  2 |\xi| +  \log  4 |\xi+\xi'| \le 2 \log |z_1| + C = 2 \log |z| + C.
\end{multline*} 

To prove that the above choice of $\zeta_j$'s allows for a map passing through the required points, we
will show that, for small enough $|\varepsilon|$, the conditions
 \eqref{ineqex1} and \eqref{ineqex2} are satisfied. 
To see this, we need to compute the quantities $w_j$, using \eqref{autos}.
\begin{eqnarray}
\label{w1} 
w_1&=& \frac{\eps_1}{\zeta_1 \phi_{\zeta_2} (\zeta_1)}=\frac{\eps_1}{\zeta_1 }\frac{1-\zeta_1 \bar \zeta_2}{\xi'}
= \frac{z_1}{\zeta_0} \frac{1-\zeta_1 \bar \zeta_2}{1-| \zeta_1 | ^2} \frac{1-\zeta_0 \bar \zeta_1}{\xi},\\
\label{w2} 
w_2&=& \frac{z_1}{\zeta_0 \phi_{\zeta_2} (\zeta_0)}=\frac{z_1}{\zeta_0} \frac{1-\zeta_0 \bar \zeta_2}{\xi+\xi'},\\
\label{w3} 
w_3&=& \frac{z_2}{\zeta_0\phi_{\zeta_1} (\zeta_0)}= \mu \frac{z_1}{\zeta_0} \frac{1-\zeta_0 \bar \zeta_1}{\xi},\\
\label{w4} 
w_4&=& \frac{\eps_2}{\zeta_2 \phi_{\zeta_1} (\zeta_2)}
=\frac{(\mu-\gamma)\eps_1}{\zeta_2 }\frac{1-\zeta_2 \bar \zeta_1}{\xi'}
\\
\nonumber
& =&
(\mu-\gamma) \frac{z_1}{\zeta_0} \frac{\zeta_1}{\zeta_2}\frac{1-\zeta_2 \bar \zeta_1}{1-| \zeta_1 | ^2} \frac{1-\zeta_0 \bar \zeta_1}{\xi}.
\end{eqnarray} 

We need to see that all those $w_j$ lie in the unit disk. It follows from the Standing assumptions above that
\begin{eqnarray}
|w_1| &\le& \left| \frac{z_1}{\xi} \right| \frac{1}{(1-| \zeta_1 | ^2)\zeta_0} \le \frac{5}{C_1} = \frac{1}{8},\\
|w_2| &\le& \left| \frac{z_1}{\xi} \right| \frac{2}{\zeta_0} \le \frac{4}{C_1} = \frac{1}{10},\\
|w_3| &\le& \left| \frac{z_1}{\xi} \right| \frac{1}{\zeta_0} \le \frac{2}{C_1} = \frac{1}{20},\\
|w_4| &\le& |\mu-\gamma| \left| \frac{z_1}{\xi} \right| \frac{2}{(1-| \zeta_1 | ^2)\zeta_0} \le \frac{20}{C_1} = \frac{1}{2},
\end{eqnarray} 
where the last inequality holds when we assume $\varepsilon$ small enough so that $|\gamma|\le 1$. It also follows from the above estimates that $|1-w_1 \bar w_2| \ge 1/2$, $|1-w_3 \bar w_4| \ge 1/2$, so that $d_G (w_1, w_2) \le 2 |w_1 - w_2|$ and $d_G (w_3, w_4) \le 2 |w_3 - w_4|$. We proceed to the computation of those Euclidean distances, which we will compare to the estimates \eqref{autoest}. 
\begin{multline*}
w_2 - w_1 = \frac{z_1}{\zeta_0} \left[ (1-\zeta_0 \bar \zeta_2) \left( \frac{1}{\xi+\xi'} -\frac{1}{\xi}\right) +\right.
\\
\left.
\frac{1}{\xi} \left( (1-\zeta_0 \bar \zeta_2)- (1-\zeta_0 \bar \zeta_1)\right) +
\frac{1-\zeta_0 \bar \zeta_1}{\xi} \left( 1- \frac{1-\zeta_1 \bar \zeta_2}{1-| \zeta_1 | ^2}\right) 
 \right] .
\end{multline*} 
Each term is estimated by 
\begin{eqnarray*}
\left| \frac{1}{\xi+\xi'} -\frac{1}{\xi} \right| &=& \left| \frac{\xi'}{\xi(\xi+\xi')}  \right| \le \frac{8|\varepsilon_1|}{C_1 |z_1|^ 2},\\
\left| (1-\zeta_0 \bar \zeta_2)- (1-\zeta_0 \bar \zeta_1)\right| &=& |\xi' \zeta_0| \le 2 |\varepsilon_1| C_1, \\
\left| 1- \frac{1-\zeta_1 \bar \zeta_2}{1-| \zeta_1 | ^2} \right| &=& \frac{|\xi' \zeta_1|}{1-| \zeta_1 | ^2} \le 2  |\xi' | 
\le 8 |\varepsilon_1| C_1,
\end{eqnarray*} 
so $|w_2 - w_1| \le 40  |\varepsilon_1| C_1 / |z_1|$. 

Similarly,
\begin{multline*}
w_3 - w_4 = \mu \frac{z_1}{\zeta_0} \frac{1-\zeta_0 \bar \zeta_1}{\xi} \left[ 1- \frac{\zeta_1}{\zeta_2} +
\frac{\zeta_1}{\zeta_2} \left(  1- \frac{1-\zeta_1 \bar \zeta_2}{1-| \zeta_1 | ^2} \right) 
 \right] \\
+
\gamma \frac{z_1}{\zeta_0} \frac{\zeta_1}{\zeta_2}  \frac{1-\zeta_1 \bar \zeta_2}{1-| \zeta_1 | ^2} \frac{1-\zeta_0 \bar \zeta_1}{\xi}.
\end{multline*} 
Using the fact that $|1- \frac{\zeta_1}{\zeta_2} |\le 8 |\xi'|$, we find 
\begin{equation*}
|w_3 - w_4| \le 96 |\varepsilon_1| + 12 |\gamma| /C_1. 
\end{equation*} 
It is then easy to see, using \eqref{autoest}, that \eqref{ineqex1} and \eqref{ineqex2} can be verified whenever we choose $\varepsilon_1$ small enough (and therefore $|\gamma|$ small enough), depending on $|z_1|$.
\end{proof*}

\bibliographystyle{amsplain}

\vskip.2cm
%\vfill\eject

Pascal J. Thomas

Laboratoire Emile Picard, UMR CNRS 5580

Universit\'e Paul Sabatier

118 Route de Narbonne

F-31062 TOULOUSE CEDEX

France

pthomas@cict.fr

\end{document}